# The existence and decay of solitary waves for the Fornberg-Whitham equation

Yong Zhang · Fei Xu · Fengquan Li

**Abstract** In this paper, we consider the Fornberg-Whitham equation and a family of solitary wave solutions is found by using minimization principle, where a related penalization function and the concentration-compactness lemma play a key role in our proof. Besides, we also prove that the family of solitary solutions is orbital stable and decay exponentially when speed wave $c$ is bigger than 1.

## 1 Introduction

The Fornberg-Whiathm(FW) equation

$$u_t - u_{txx} + u_x + uu_x = uu_{xxx} + 3u_x u_{xx} \tag{1.1}$$

was derived by B.Fornberg and G.B. Whitham as a model to study the qualitative behaviors of wave-breaking in [1], where they also found a special peaked solitary wave solution $u(t,x) = Ae^{-\frac{1}{2}|x-\frac{4}{3}t|}$. Since that time, a little attention was paid to it. Until the last decade, the authors in [2] constructed some new types of travelling wave solutions for (1.1) and made a classification of its travelling waves in [3]. There are gradually growing appeals for FW equation. In [4][5], several blow-up phenomena of the Fornberg-whitham equation on line $R$ and on circle $T$ were established. The author in [6] employed a Galerkin type approximation

―――――――――――――――
Y. Zhang · F. Xu · F. Li(✉)
School of Mathematical Sciences, Dalian University of Technology, Dalian, 116024, China
e-mail: fqli@dlut.edu.cn

argument showing that its Cauchy problem was well-posed in Sobolev spaces $H^s(T)$, $s > \frac{3}{2}$. The existence of global attractor to the viscous Fornberg-Whitham equation was proved in [7] and discontinuous traveling waves as weak solutions to the Fornberg-Whitham equation were investigated in [8]. Note that some new results on two-component Fornberg-Whitham equation can be seen in [9]. Despite many endeavors, some important properties of (1.1) are still unclear.

Unlike Camassa-Holm (CH) equation or Degasperis-Procesi (DP) equation, despite they seem to only differ in coefficients, FW equation loses the property of complete integrability and some good conversation laws (see [10]). On the other hand, if we write FW equation in following nonlocal form

$$u_t + uu_x + (1 - \partial_x^2)^{-1} u_x = 0 \tag{1.2}$$

to investigate its travelling wave solutions. We would find that there is a great difference between KdV-type equations in [11] (or Benjamin-type in [12] and BBM-type equations in [13]) and (1.2), where the operator $L := (1 - \partial_x^2)^{-1}$ is not of positive order. These special characteristics of FW equation add a lot of difficulties to our study, which compels us to find new approaches.

The main contributions of our work can be summarised as follows. In Section 3, we introduce the penalization function and then use the constrained minimization principle to prove the existence of periodic waves of FW equation. Based on these periodic minimizers, we construct a global minimizing sequence and find the global minimizer (i.e solitary wave) by using the concentration-compactness lemma in Section 4. The last Section is dedicated to prove that these solitary waves are orbital stable and have a exponential decay.

## 2 Preliminaries and main results

Note that FW equation can be expressed in following nonlocal form

$$u_t + uu_x + Lu_x = 0, \tag{2.1}$$

where $\widehat{Lf}(\xi) = \frac{1}{1+\xi^2} \hat{f}(\xi)$. Since considering the travelling wave solutions to (2.1), we may assume that

$$u(t, x) = u(x - ct), \tag{2.2}$$

where $c$ is wave speed. And let $y = x - ct$, then (2.1) and (2.2) imply that

$$-cu + \frac{u^2}{2} + Lu = 0. \tag{2.3}$$

Now define the functional

$$J(u) := -\frac{1}{2}\int_R uLu\,dy - \frac{1}{6}\int_R u^3 dy \tag{2.4}$$

and the $L^2$ constraint

$$Q(u) := \int_R u^2 dy = q > 0. \tag{2.5}$$

**Remark 2.1.** *If we can find a minimizer $\bar{u} \in H^1(R)$ of the functional $J(u)$ over the constraint space $\{u \in H^1(R)|Q(u) = q\}$, then the Lagrange multiplier theorem implies that there exists a $\lambda$ such that*

$$J'(\bar{u})v = \lambda Q'(\bar{u})v, \quad for\ v \in C_c^\infty(R).$$

*That is to say,*

$$-\int_R L\bar{u}v\,dy - \frac{1}{2}\int_R \bar{u}^2 v\,dy = 2\lambda \int_R \bar{u}v\,dy.$$

*If we choose the Lagrange multiplier $\lambda = -\frac{c}{2}$, then it's obvious that $\bar{u} \in H^1(R)$ is a weak solution to (2.3). And we let the set $G_q$ denote the set of minimizers of $J(u)$ over corresponding constrained space. However, it's not easy to find directly this minimizer on the whole space. There are two essential difficulties to overcome, one of which is to compensate the coerciveness of $J(u)$ in $H^1(R)$ and the other is to recover the compactness on R.*

Next we shall state the main results of this paper. The following Theorem 2.1 is on the existence of solitary waves for FW equation whose proof can be found in Section 3 and Section 4. Theorem 2.2 is dedicated to the orbital stability and decay, which have been proved in Section 5.

**Theorem 2.1.** *The set $G_q$ is not empty, and each element of $G_q$ is a solitary-wave solution for equation (1.1) and $-\frac{c}{2}$ is corresponding Lagrange multiplier.*

**Theorem 2.2.** *$G_q$ is a stable set in the following sense: for every $\varepsilon > 0$ and $g \in G_q$, there exists a $\delta > 0$ such that $\|u_0 - g\|_{H^{1-\varepsilon}} < \delta$, then the solution $u \in C((0, T); H^1(R))$ of (1.1) with $u(0, x) = u_0(x)$ satisfies*

$$\inf_{g \in G_q} \|u(t, \cdot) - g\|_{H^{1-\varepsilon}} < \varepsilon, \quad for\ all\ t \in (0, T).$$

*And if the speed wave c > 1, there exists some constant C ∈ R such that*

$$\lim_{y\to\pm\infty} e^{\sqrt{\frac{c-1}{c}}|y|} g(y) = C.$$

# 3 Existence of periodic waves

In this section, the method developed in [14] and [15] is not suitable due to the loss of positive order for operator $L$. Inspired by [16, 17], here we introduce the penalization function and then use the constrained minimization principle to prove the existence of periodic waves of FW equation.

Firstly, we study a related penalization function acting on periodic functions. Let $L_P^2(-P, P]$ be the space of 2P-periodic, locally square integrable functions with fourier series representation

$$w(y) = \frac{1}{\sqrt{2P}} \sum_{k\in Z} \hat{w}(\xi) e^{\frac{i\pi y\xi}{P}}$$

where

$$\hat{w}(\xi) = \frac{1}{\sqrt{2P}} \int_{-P}^{P} w(y) e^{\frac{-i\pi y\xi}{P}} dy.$$

And we define

$$H_P^1(-P, P] = \{w \in L_P^2(-P, P] : \|w\|_{H_P^1} = \left[\sum_{k\in Z}(1 + \frac{\pi^2\xi^2}{P^2})|\hat{w}(\xi)|^2\right]^{\frac{1}{2}} < \infty\}$$

Now we consider following penalized periodic problem over a constraint space. Let the modified functional be

$$J(u|P, \rho) := -\frac{1}{2}\int_{-P}^{P} uLu\,dy - \frac{1}{6}\int_{-P}^{P} u^3 dy + \rho(\|u\|_{H_P^1}^2), \tag{3.1}$$

where $\rho : [0, 2R^2) \to [0, +\infty)$ is smooth and increasing penalization function with $\rho(t) \to 0$ as $t \in [0, R^2]$ and $\rho(t) \to +\infty$ as $t \to 2R^2$. And define the constraint space

$$V_{P,q,\sqrt{2}r} := \{u \in H_P^1 \mid Q(u|P) = \int_{-P}^{P} u^2 dy = q, \ \|u\|_{H_P^1} < \sqrt{2}r\}. \tag{3.2}$$

The existence of minimizer of $J(u|P, \rho)$ over $V_{P,q,\sqrt{2}r}$ is obvious. Because it's easy to check the new functional $J(u|P, \rho)$ is coercive and the fact that we are now working in space $H_P^1(-P, P]$ allows the use of the Rellich-Kondrachov theorem. Thus there exists $u^* \in H_P^1(-P, P]$ such that

$$J(u^*|P, \rho) = \inf_{u \in V_{P,q,\sqrt{2}r}} J(u|P, \rho). \tag{3.3}$$

Besides, here we claim that $u^*(y)$ found in (3.3) also minimizes

$$J(u|P) := -\frac{1}{2}\int_{-P}^{P} uLu \, dy - \frac{1}{6}\int_{-P}^{P} u^3 \, dy \tag{3.4}$$

over

$$V_{P,q,r} := \{u \in H_P^1 \mid Q(u|P) = \int_{-P}^{P} u^2 \, dy = q, \ \|u\|_{H_P^1} < r\}. \tag{3.5}$$

The definition of the penalization function $\rho(t)$ and (3.3) show that it's immediate if we can prove

$$\|u^*\|_{H_P^1}^2 \leq C_q. \tag{3.6}$$

Based on (3.3), we apply the Lagrange multiplier theorem to obtain

$$J'(u^*|P,\rho)v = -\frac{c}{2}Q'(u^*|P)v, \quad for \ v \in C_c^\infty(-P,P].$$

That is to say,

$$-\int_{-P}^{P}(Lu^* + \frac{(u^*)^2}{2})v \, dy + 2\rho'(\|u^*\|_{H_P^1}^2)\int_{-P}^{P} u_y^* v_y + u^* v \, dy = -c\int_{-P}^{P} u^* v \, dy.$$

Thus the following equality holds in distribution sense

$$Lu^* + \frac{(u^*)^2}{2} - 2\rho'(\|u^*\|_{H_P^1}^2)(u^* - u_{yy}^*) = cu^*. \tag{3.7}$$

Multiplying (3.7) by $u^* - u_{yy}^*$ and integrating on $(-P, P]$, we have

$$\begin{aligned}
c\|u^*\|_{H_P^1}^2 &= \int_{-P}^{P}(Lu^* + \frac{(u^*)^2}{2})(u^* - u_{yy}^*)dy - 2\rho'(\|u^*\|_{H_P^1}^2)\int_{-P}^{P}(u^* - u_{yy}^*)^2 dy \\
&\leq \int_{-P}^{P}(Lu^* + \frac{(u^*)^2}{2})(u^* - u_{yy}^*)dy \\
&= \int_{-P}^{P} u^* Lu^* + u_y^* Lu_y^* dy + \int_{-P}^{P} \frac{(u^*)^3}{2} + u^*(u_y^*)^2 dy \\
&\leq 2\|u^*\|_{L_P^2}^2 + \|u^*\|_{L_P^\infty}\|u^*\|_{H_P^1}^2,
\end{aligned} \tag{3.8}$$

where the last inequality uses the Plancherel equality and the symbol of operator $L$ is less than one. Because $\|u^*\|_{L_P^\infty} \to 0$ as $\|u^*\|_{H_P^1} \to 0$, then $\|u^*\|_{L_P^\infty}$ is bounded by $\frac{c}{2}$ for small values of $r$. Thus this implies that

$$\|u^*\|_{H_P^1}^2 \leq \frac{4}{c}\|u^*\|_{L_P^2}^2 = \frac{4}{c}q$$

Up to now, we have proved the existence of periodic waves to (2.3).

# 4 Existence of solitary waves

In this section, we would prove that

$$I_q := \inf_{u \in V_{q,r}} J(u)$$

can be attained at some $\bar{u} \in V_{q,r}$. And one of key ingredients is to construct a minimizing sequence for $J(u)$. The main approach is to extend the minimizer $u^*$ of $J(u|P)$ over $V_{P,q,r}$ by translation and truncation in following manner.

Let $\{u_n^*\}_n$ be a bounded family of functions in $H^1(R)$ with

$$Supp\ \{u_n^*\} \subset (-P, P)\ \ and\ \ u_n^* \to u^*\ in\ L^2(R)\ as\ n \to \infty. \tag{4.1}$$

Define

$$\bar{u}_n = \sum_{j \in Z} u_n^*(\cdot + 2Pj), \tag{4.2}$$

now we turn to prove that there exists a subsequence of $\{u_{n_k}^*\}_k$, denote $\{u_n^*\}_n$ again, such that

$$J(\bar{u}_n) = \lim_{P \to \infty} J(u_n^*|P),\quad Q(\bar{u}_n) = \lim_{P \to \infty} Q(u_n^*|P). \tag{4.3}$$

Firstly we observe that

$$\begin{aligned}\|L\bar{u}_n - Lu_n^*\|_{L_P^2}^2 &= \int_{-P}^{P} |L\bar{u}_n - Lu_n^*|^2 dy \\ &\leq \int_{-P}^{P} |\bar{u}_n - u_n^*|^2 dy \\ &= \int_{-P}^{P} |\sum_{|j| \geq 1} u_n^*(y + 2Pj)|^2 dy = 0,\end{aligned} \tag{4.4}$$

due to $Supp\ \{u_n^*(y + 2Pj)\} \subset (-P + 2Pj, P + 2Pj)$. Similarly, we also have

$$\|\bar{u}_n - u_n^*\|_{L_P^2}^2 = 0. \tag{4.5}$$

Then we can prove that

$$\begin{aligned}&|\frac{1}{2}\int_{-P}^{P} \bar{u}_n L\bar{u}_n - u_n^* Lu_n^* dy + \frac{1}{6}\int_{-P}^{P} (\bar{u}_n)^3 - (u_n^*)^3 dy| \\ &\leq \frac{1}{2}\int_{-P}^{P} |\bar{u}_n L\bar{u}_n - \bar{u}_n Lu_n^* + \bar{u}_n Lu_n^* - u_n^* Lu_n^*| dy + \frac{1}{6}\int_{-P}^{P} |(\bar{u}_n - u_n^*)[(\bar{u}_n)^2 + (u_n^*)^2 + \bar{u}_n u_n^*]| dy \\ &\leq \frac{1}{2}(\|\bar{u}_n\|_{L_P^2}\|L\bar{u}_n - Lu_n^*\|_{L_P^2} + \|Lu_n^*\|_{L_P^2}\|\bar{u}_n - u_n^*\|_{L_P^2}) + \frac{1}{6}\|\bar{u}_n - u_n^*\|_{L_P^2}\|(\bar{u}_n)^2 + (u_n^*)^2 + \bar{u}_n u_n^*\|_{L_P^2} \\ &\to 0,\end{aligned} \tag{4.6}$$

as $P \to \infty$. The last equality uses (4.4), (4.5) and the boundedness of $\{\bar{u}_n\}$ and $\{u_n^*\}$. Thus we finish the proof of (4.3), and we can conclude that $\{\bar{u}_n\}$ is the minimizing sequence for $J(u)$ over $V_{q,r}$.

The other key is to recover the compactness by using the concentration-compactness lemma in [18]. For convenience, we here state it for our purpose.

**Lemma 4.1.** *Let $\{\eta_n\}_{n\geq 1}$ be a sequence of nonnegative function in $L^1(R)$ with $\int_R \eta_n(x)dx = q$ for all n and some $q > 0$. Then there exists a subsequence $\{\eta_{n_k}\}_{k\geq 1}$, denote $\{\eta_n\}_{n\geq 1}$ again, satisfying one of the following cases.*

*Case1: Vanashing*

*For each $r > 0$, then*

$$\lim_{n \to \infty} sup_{x \in R} \int_{B_r(x)} \eta_n(y)dy = 0.$$

*Case2: Concentration*

*There is a sequence $\{x_n\}_{n \in N} \subset R$, such that for $\varepsilon > 0$ and $\exists r' > 0$*

$$\int_{B_{r'}(x_n)} \eta_n(y)dy \geq q - \varepsilon.$$

*Case3: Dichotomy*

*There exists $p \in (0, q)$, $M \geq 1$ and nonnegative $\eta_n^{(1)}, \eta_n^{(2)} \in L^1(R)$, s.t. for any $\varepsilon > 0, n \geq M$*

$$\begin{cases} \|\eta_n - (\eta_n^{(1)} + \eta_n^{(2)})\|_{L^1} \leq \varepsilon, \\ |\int_R \eta_n^{(1)} dy - p| \leq \varepsilon, \quad |\int_R \eta_n^{(2)} dy - (q-p)| \leq \varepsilon, \\ Supp\, \eta_n^{(1)} \cap Supp\, \eta_n^{(2)} = \emptyset, \quad dist\{Supp\, \eta_n^{(1)}, Supp\, \eta_n^{(2)}\} \to \infty \text{ as } n \to \infty. \end{cases}$$

Now let $\eta_n = \bar{u}_n^2$, $n \in N$, where $\bar{u}_n$ is a minimizing sequence for $J(u)$ over $V_{q,r}$. Firstly it's straightforward to exclude the Case1. Indeed, if Case1 occurs for any $r' > 0$, then

$$0 < q = Q(\bar{u}_n) = \frac{1}{2} \int_R \bar{u}_n^2 dy \leq \lim_{n \to \infty} sup_{x \in R} \int_{B_{r'}(x)} \eta_n(y)dy = 0,$$

which is a contradiction.

If the Case2 occurs, then from

$$\bar{u}_n \rightharpoonup \bar{u} \text{ in } H^1(R),$$

we can infer
$$\bar{u}_n \to \bar{u} \quad in \quad L^2(R).$$
That is to say, $\bar{u} \in V_{q,r}$ and $I_q = J(\bar{u})$, then we finish the proof.

However, it's tricky to exclude the Case3. For a clear presentation, let us divide this process into the following three steps.

**Step1:** If Dichotomy occurs, we would find $p \in (0, q)$, s.t
$$I_q \geq I_p + I_{q-p}. \tag{4.7}$$

pf: Without loss of generality, we may assume that $\eta_{n_k}^{(1)}$ and $\eta_{n_k}^{(2)}$ satisfy
$$Supp\ \eta_n^{(1)} \subset (x_n - R_n, x_n + R_n), \quad Supp\ \eta_n^{(2)} \subset (-\infty, x_n - 4R_n) \cup (x_n + 4R_n, \infty),$$
where $x_n \in R$ and $R_n \to \infty$. Then it's obvious that
$$\int_{R_n \leq |y-x_n| \leq 4R_n} \eta_n dy \leq \varepsilon,$$
i.e.
$$\int_{R_n \leq |y-x_n| \leq 4R_n} \bar{u}_n^2 dy \leq \varepsilon. \tag{4.8}$$

Choose $\varphi, \psi, \phi \in C^\infty(R)$ such that $0 \leq \varphi, \psi, \phi \leq 1$ for all $y$, with
$$\varphi(y) = \begin{cases} 1, & |y| \leq 2 \\ 0, & |y| \geq 3 \end{cases}, \quad \psi(y) = \begin{cases} 1, & |y| \geq 3 \\ 0, & |y| \leq 2 \end{cases}, \quad \phi(y) = \begin{cases} 1, & 2 \leq |y| \leq 3 \\ 0, & |y| \leq 1 \ and\ |y| \geq 4. \end{cases}$$

Let $\varphi_n(y) = \varphi(\frac{y-x_n}{R_n}), \psi_n(y) = \psi(\frac{y-x_n}{R_n})$ and $\phi_n(y) = \phi(\frac{y-x_n}{R_n})$, define
$$\bar{u}_n^{(1)} = \varphi_n \bar{u}_n, \quad \bar{u}_n^{(2)} = \psi_n \bar{u}_n, \quad w_n = \phi_n \bar{u}_n$$

Since $Q(\bar{u}_n^{(1)})$ is bounded, there is a subsequence of $\bar{u}_n^{(1)}$, still denoted by $\{\bar{u}_n^{(1)}\}$, and a $p$ such that
$$Q(\bar{u}_n^{(1)}) = p, \quad for\ large\ enough\ n. \tag{4.9}$$

On the other hand, we have
$$q = Q(\bar{u}_n) = \int_R \bar{u}_n^2 dy = \int_{|y-x_n| \leq 2R_n} \bar{u}_n^2 dy + \int_{|y-x_n| \geq 3R_n} \bar{u}_n^2 dy + \int_{2R_n \leq |y-x_n| \leq 3R_n} \bar{u}_n^2 dy$$
$$= \int_{|y-x_n| \leq 2R_n} (\bar{u}_n^{(1)})^2 dy + \int_{|y-x_n| \geq 3R_n} (\bar{u}_n^{(2)})^2 dy + \int_{2R_n \leq |y-x_n| \leq 3R_n} \bar{u}_n^2 dy$$
$$= Q(\bar{u}_n^{(1)}) + Q(\bar{u}_n^{(2)}) + O(\varepsilon). \tag{4.10}$$

The last equality uses (4.8). And (4.9) and (4.10) show that

$$Q(\bar{u}_n^{(2)}) = q - p, \ for \ large \ enough \ n. \tag{4.11}$$

To finish this claim, we observe that

$$\begin{aligned} J(\bar{u}_n) &= J(\bar{u}_n^{(1)} + \bar{u}_n^{(2)} + w_n) \\ &= -\frac{1}{2}\int_R (\bar{u}_n^{(1)} + \bar{u}_n^{(2)} + w_n)L(\bar{u}_n^{(1)} + \bar{u}_n^{(2)} + w_n)dy - \frac{1}{6}\int_R (\bar{u}_n^{(1)} + \bar{u}_n^{(2)} + w_n)^3 dy \\ &= J(\bar{u}_n^{(1)}) + J(\bar{u}_n^{(2)}). \end{aligned} \tag{4.12}$$

The last equality uses adequately the definition of $\bar{u}_n^{(1)}, \bar{u}_n^{(2)}$ and $w_n$. Thus (4.9), (4.11) and (4.12) imply that

$$I_q = \lim_{n\to\infty} J(\bar{u}_n) \geq \lim_{n\to\infty} J(\bar{u}_n^{(1)}) + \lim_{n\to\infty} J(\bar{u}_n^{(2)}) = I_p + I_{q-p}.$$

**Step2:** Here we assert that

$$\int_R (\bar{u}_n)^3 dy > 0. \tag{4.13}$$

pf: Here we prove this assertion by method of contradiction, i.e. assume that

$$\int_R (\bar{u}_n)^3 dy \leq 0. \tag{4.14}$$

Firstly because $\{\bar{u}_n\}_n$ is the minimizing sequence of $J(u)$ over $V_{q,r}$, then it's also minimizing sequence for

$$\bar{I}_q := \inf_{V_{q,r}} \bar{J}(u),$$

where $\bar{J}(u) = J(u) + q = J(u) + Q(u) = \int_R u(1 - \frac{L}{2})u dy - \frac{1}{6}\int_R u^3 dy$. Then

$$\bar{I}_q = \lim_{n\to\infty} \bar{J}(\bar{u}_n) = \int_R \bar{u}_n(1 - \frac{L}{2})\bar{u}_n dy - \frac{1}{6}\int_R (\bar{u}_n)^3 dy \geq -\frac{1}{6}\int_R (\bar{u}_n)^3 dy \geq 0, \tag{4.15}$$

where the first inequality uses the symbol of $L$ is less than 1 and the last uses (4.14).

On the other hand, we can choose a function $h \in H^1(R)$ such that $Q(h) = \int_R h^2 dx = q$, $\|h\|_{H^1} < r$ and $\int_R h^3 dy > 0$. Then for any $0 < \theta < 1$, we have

$$Q(\sqrt{\theta}h(\theta y)) = q, \quad \|\sqrt{\theta}h(\theta y)\|_{H^1} < r.$$

That is to say, $\sqrt{\theta}h(\theta y) \in V_{q,r}$. At the same time,

$$\begin{aligned} \bar{J}(\sqrt{\theta}h(\theta y)) &= \int_R \sqrt{\theta}h(\theta y)(1 - \frac{L}{2})\sqrt{\theta}h(\theta y)dy - \frac{1}{6}\int_R (\sqrt{\theta}h(\theta y))^3 dy \\ &= \int_R h(y)(1 - \frac{L}{2})h(y)dy - \frac{\sqrt{\theta}}{6}\int_R (h(y))^3 dy < 0, \end{aligned} \tag{4.16}$$

for suitable choice of $q$ and $\theta$. This indicates that $\bar{I}_q < 0$, which is contradicted with (4.15). Thus we have $\int_R (\bar{u}_n)^3 dy > 0$.

**Step3:** We claim that $I_q < I_p + I_{q-p}$.

Because $I_q = \inf_{V_{q,r}} J(\bar{u}_n)$, then

$$\gamma I_q = \gamma \inf_{V_{q,r}} J(\bar{u}_n) = -\frac{\gamma}{2} \int_R \bar{u}_n L \bar{u}_n dy - \frac{\gamma}{6} \int_R (\bar{u}_n)^3 dy, \tag{4.17}$$

and

$$I_{\gamma q} = \inf_{V_{q,r}} J(\sqrt{\gamma} \bar{u}_n) = -\frac{\gamma}{2} \int_R \bar{u}_n L \bar{u}_n dy - \frac{\gamma \sqrt{\gamma}}{6} \int_R (\bar{u}_n)^3 dy. \tag{4.18}$$

Combining (4.17) (4.18) with (4.13), we can obtain for any $\gamma > 1$

$$I_{\gamma q} < \gamma I_q. \tag{4.19}$$

Without loss of generality, we assume that $p > q - p$, then (4.19) indicates that

$$I_q = I_{p(1+\frac{q-p}{p})} < (1 + \frac{q-p}{p})I_p = I_p + \frac{q-p}{p}I_p = I_p + \frac{q-p}{p}I_{\frac{p}{q-p}(q-p)} < I_p + I_{q-p}.$$

This is contradicted with (4.7). Up to now, we exclude the Case3 and finish the proof of solitary waves' existence.

# 5 Orbital stability and decay of solitary waves

In this section, we are devoted to orbital stability for solitary waves found in Section 4 and use Bona and Li's lemma in [19] to infer that the solitary waves decay at infinite points with exponential form if the wave speed $c > 1$.

According to the analysis in Section 4, we know that if $\{v_n\}_n$ is a minimizing sequence of $J(u)$ and it is uniformly bounded in $H^1(R)$, then $\{v_n\}_n$ converges strongly in $L^2(R)$. In fact, it strongly converges in all the intermediate norms according to interpolation lemma in Sobolev spaces. That is to say, there is a sequence of real numbers $\{z_n\}$, such that $\{v_n(\cdot + z_n)\}$ has a subsequence that converges in $H^{1-\varepsilon}(R)$ to an element $g \in G_q$

$$\|v_n(\cdot + z_n) - g\|_{H^{1-\varepsilon}} \le C_\varepsilon \|v_n(\cdot + z_n)_n - g\|_{L^2}^\varepsilon \|v_n(\cdot + z_n) - g\|_{H^1}^{1-\varepsilon} \to 0, \quad as\ n \to \infty. \tag{5.1}$$

Since $g(\cdot + z)$ is also in $G_q$, a simple translation and (5.1) imply that

$$\lim_{n \to \infty} \inf_{g \in G_q} \|v_n - g\|_{H^{1-\varepsilon}} = 0 \tag{5.2}$$

Assume that the Theorem 2.2 is false, then there exists a $g^* \in G_q$ and $\varepsilon > 0$, such that for every $n \in N$, we can find $\phi_n \in H^{1-\varepsilon}(R)$ and $t_n \in R$ satisfy

$$\|\phi_n - g^*\|_{H^{1-\varepsilon}} \leq \frac{1}{n}, \tag{5.3}$$

and

$$\inf_{g \in G_q} \|\bar{u}_n(t_n, \cdot) - g\|_{H^{1-\varepsilon}} \geq \varepsilon, \tag{5.4}$$

where $\bar{u}_n(t_n, \cdot)$ is the solution of (1.1) with $\bar{u}_n(0, \cdot) = \phi_n$. Then (5.3) and Sobolev embedding theorem imply that

$$J(\phi_n) \to J(g^*) = I_q, \quad q_n := Q(\phi_n) \to Q(g^*) = q.$$

Define $v_n := \sqrt{\frac{q}{q_n}} \bar{u}_n(t_n, \cdot)$, then

$$Q(v_n) = \frac{q}{q_n} Q(\bar{u}_n(t_n, \cdot)) = \frac{q}{q_n} Q(\phi_n) \to q, \tag{5.5}$$

and

$$|J(v_n) - J(\phi_n)| = |J(v_n) - J(\bar{u}_n(t_n, \cdot))|$$
$$= |\frac{1}{2}\left(1 - \frac{q}{q_n}\right) \int_R \bar{u}_n L \bar{u}_n dx + \frac{1}{6}\left[1 - \left(\frac{q}{q_n}\right)^{\frac{3}{2}}\right] \int_R \bar{u}_n^3 dx|$$
$$\leq C \|\bar{u}_n(t_n, \cdot)\|_{H^1}(|1 - \frac{q}{q_n}| + |1 - \left(\frac{q}{q_n}\right)^{\frac{3}{2}}|)$$
$$\leq Cr(|1 - \frac{q}{q_n}| + |1 - \left(\frac{q}{q_n}\right)^{\frac{3}{2}}|) \to 0, \tag{5.6}$$

as $n \to \infty$. Thus the sequence $\{v_n\}_n$ is also a minimizing sequence of $J(u)$ over $V_{q,r}$. Therefore, we can infer that

$$\varepsilon \leq \inf_{g \in G_q} \|\bar{u}_n(t_n, \cdot) - g\|_{H^{1-\varepsilon}} \leq \|\sqrt{\frac{q}{q_n}} \bar{u}_n(t_n, \cdot) - g^*\|_{H^{1-\varepsilon}} = \|v_n - g^*\|_{H^{1-\varepsilon}} = \inf_{g \in G_q} \|v_n(t_n, \cdot) - g\|_{H^{1-\varepsilon}} \leq \frac{\varepsilon}{2},$$

where the first inequality uses (5.4), the second equality uses (5.5) and (5.6), the last inequality uses (5.2). This leads to a contradiction. Up to now, we have proved the orbital stability.

At last, we introduce the Bona and Li's lemma to show the asymptotic decay of this solitary-wave solutions.

**Lemma 5.1.** *Suppose that $f \in L^{\infty}(R)$ with $\lim_{y \to \pm\infty} f(y) = 0$ is a solution of the convolution equation*

$$f(y) = \int_R k(y-x)M(f(x))dx,$$

*where the kernel $k$ is a measurable function satisfying $\hat{f} \in H^s$, $s > \frac{1}{2}$ and $M$ is a function such that $|M(t)| \leq C|t|^m$ for all $t \in R$ for some constants $C > 0$ and $m > 1$. Then*

*(1) $f \in L^1(R)$;*

*(2) there is a constant $l$ with $0 < l < s$ such that $|y|^l f(y) \in L^1 \cap L^{\infty}$;*

*(3) if $\lim_{y \to +\infty} e^{\sigma y} k(y) = C^*$ for some constants $C^*$ and $\sigma > 0$, then*

$$\sup_{y \in R} e^{\sigma |y|} f(y) = C.$$

Based on the travelling form (2.3), we can express FW equation in following convolution form

$$u = \frac{1}{2}(c-L)^{-1}u^2, \tag{5.7}$$

where the kernel $k(x) = F^{-1}(\frac{\xi^2+1}{c\xi^2+c-1})$ and $M(u) = \frac{1}{2}u^2$. According to Lemma 5.1, it's sufficient for us to prove

$$\lim_{y \to +\infty} e^{\sigma y} k(y) = C^*, \quad for \ some \ \sigma > 0.$$

In fact, if $c > 1$, we have

$$\begin{aligned}
k(y) &= \frac{1}{\sqrt{2\pi}} \int_R \frac{\xi^2+1}{c\xi^2+c-1} e^{i\xi y} d\xi \\
&= \frac{1}{\sqrt{2\pi}c} \check{1} + \frac{1}{\sqrt{2\pi}} \int_R \frac{e^{i\xi y}}{c^2\xi^2+c(c-1)} d\xi \\
&= \frac{1}{\sqrt{2\pi}c} \delta(y) + \sqrt{2\pi} i \ Res \ [\frac{e^{i\xi y}}{c^2\xi^2+c(c-1)}, i\sqrt{\frac{c-1}{c}}] \\
&= \frac{1}{\sqrt{2\pi}c} \delta(y) + \frac{\sqrt{2\pi}}{2\sqrt{c^3(c-1)}} e^{-\sqrt{\frac{c-1}{c}}y} \\
&= \frac{\sqrt{2\pi}}{2\sqrt{c^3(c-1)}} e^{-\sqrt{\frac{c-1}{c}}y}, \quad for \ y \neq 0
\end{aligned}$$

where $\delta(y)$ is the dirac function and we use Jordan's lemma and Residue theorem in complex analysis. Up to now, we finish the proof if choose $\sigma = \sqrt{\frac{c-1}{c}}$ and $C^* = \frac{\sqrt{2\pi}}{2\sqrt{c^3(c-1)}}$.

**Remark 5.1.** *In fact, the special solitary wave $u(t,x) = Ae^{-\frac{1}{2}|x-\frac{4}{3}t|}$ found by Fornberg and Whitham in [1] belongs to $G_q$. At this time, we only need to let $c = \frac{4}{3}$, then $\sigma = \sqrt{\frac{c-1}{c}} = \frac{1}{2}$.*

*However, we can't infer that if the solitary waves, whose wave speed is not greater than one, exist. And if there exist such solitary waves, then how about its rate of decay? In our future work, we plan to investigate these interesting problems.*

# Acknowledgements

This project is supported by National Natural Science Foundation of China (No:11571057).